\documentclass[10pt]{article}
\usepackage[latin1]{inputenc}
\usepackage[T1]{fontenc}
\usepackage[english]{babel}

\usepackage{amssymb, amsmath, amsthm ,amsxtra,amsfonts, euscript,latexsym, dsfont, mathrsfs,
eufrak, yfonts, longtable,  textcomp, marvosym}
\usepackage{graphics,graphicx}
\usepackage{epsf}
\usepackage{enumerate}
\usepackage{array}

\renewcommand{\sectionmark}[1]%
{\markright{\MakeUppercase{\thesection.\#1}}}

\theoremstyle{Theorem}
\numberwithin{equation}{section}

\newtheorem{teo}{Theorem}[section]
\newtheorem*{teo1}{Theorem}

\newtheorem{defi}[teo]{Definition}
\newtheorem{os}[teo]{Remark}
\newtheorem{lemma}[teo]{Lemma}

\newcommand{\rn}{$\mathbb{R}^N$}
\newcommand{\front}{\partial{\Omega}}

\newcommand{\dn}{\displaystyle{\frac{\partial{u}}{\partial{n}}}}
\newcommand{\ins}{$\Omega$}
\newcommand{\cont}{$C^{k}(\overline{\Omega})$}

\newcommand{\hold}{$C^{k,\alpha}(\overline{\Omega})$}

\newcommand{\eqz}{C^{0,\alpha}(\overline{\Omega})}

\newcommand{\borneu}{C^{1,\alpha}(\overline{\Omega})}
\newcommand{\sol}{C^{2,\alpha}(\overline{\Omega})}
\newcommand{\dimo}{{\bf Proof}}

\renewcommand{\sectionmark}[1]%
{\markright{\MakeUppercase{\thesection.\#1}}}

\title{Schauder estimate for solutions of Poisson's equation with  Neumann boundary condition}

\author{
G. Nardi\footnote{CEntre de REcherches en MAth\'ematiques de la DEcision (CEREMADE), CNRS : UMR7534,  Universit\'e Paris IX - Paris Dauphine.  Email: {\tt nardi@ceremade.dauphine.fr}} 
}

\date{}

\begin{document}

\maketitle

\begin{abstract}
In this work we consider the Neumann problem for the Laplace operator and we prove an existence result in the H\"older spaces and obtain Schauder estimates. According to our knowledge this result is not explicitly  proved in the several works devoted to the Schauder  theory, where similar theorems are proved for the Dirichlet and oblique derivative problems. Our contribution is to make explicit the existence and the estimate for the Neumann problem.

\end{abstract}

\section{Introduction}

Let $\Omega$ be an $C^{2,\alpha}$-domain of \rn (we refer to Section \ref{GD} for notation and definitions). We consider the following problem 
$$
\left\{
\begin{array}{rl}
\Delta u=f &\textrm{ in $\Omega$} \\
\dn=g &\textrm{ on $\front$}
\end{array} \right.
$$
with $f\in C^{0,\alpha}(\overline{ \Omega})$ and  $g\in\borneu$. The aim of this paper is to prove an existence result in $\sol$ for this problem and an  estimate of the form
$$\|u\|_{C^{2,\alpha}}\leq C(\|f\|_{C^{0,\alpha}}+\|g\|_{C^{1,\alpha}}).$$
In the 1930s, this kind of estimates has been used by Schauder \cite{S1} and Caccioppoli \cite{C} to prove an existence result in $\sol$ for the Dirichlet problem for an elliptic equation (\cite{GT}: Theorem 6.8 p.100; \cite{LU}: Theorem 1.3 p.107). Using the same technique,  in the 1950s, Fiorenza \cite{F} proved a similar estimate and an existence result in $\sol$ for the oblique derivative problem 
$$l(x)u+m(x)\dn=g \;\;\mbox{on}\;\;\front \quad\quad (lm>0 \;\mbox{on}\;\front )$$ 
for elliptic equations (\cite{GT}:  Theorem 6.31, p.128; \cite{LU}: Theorem 3.1, p.126). 
 \par Unfortunately the hypothesis $l\neq 0$ can not be removed in the proof of the existence result and the Schauder estimate  for the oblique derivative problem (\cite{GT}: Theorem 6.31, p.128). Moreover, reading this proof, we can verify that it is not even possible to get the result for the Neumann problem taking the limit $l\rightarrow 0$.  
 \par So the case of the Neumann problem needs to be considered independently and, according to  our knowledge, it is not explicitly present in the classical litterature on the subject (see for instance \cite{GT}, \cite{LU}) where the Dirichlet and oblique derivative problem are studied in detail.
\par  The main goal of this work is to  formalize the existence result and a Schauder estimate for the Neumann problem for the Poisson's equation.
\par Our main contribution is the following result (Theorems \ref{esistenza} and \ref{stima}):
\begin{teo1} Let  $\alpha\in (0,1)$  and   $\Omega$ be a  $C^{2,\alpha}$-domain. Let  $f\in \eqz$ and   $g\in\borneu$ such that
$$\int_{\Omega}f=\int_{\partial \Omega}g.$$
Then there exists  a solution $u\in \sol$  (unique up to an additive constant) to the problem
$$
\left\{
\begin{array}{rl}
\Delta u=f &\textrm{ in $\Omega$} \\
\dn=g &\textrm{ on $\front$}
\end{array} \right.
$$
Moreover, every solution to this problem verifies the following estimate
$$\left\|u-\frac{1}{|\Omega|}\int_{\Omega} u\right\|_{C^{2,\alpha}}\leq C(\|f\|_{C^{0,\alpha}}+\|g\|_{C^{1,\alpha}}),$$
with $C=C(\Omega,N,\alpha).$
\end{teo1}
The starting point of the proof is the alternative theorem for the oblique derivative problem  for an uniformly elliptic operator $L$ with coefficients in  $\eqz$ and   $c< 0$ (Theorem \ref{alternativa}). This allows us to prove an existence (uniqueness is given up to a constant) result for the Neumann problem for elliptic operators with $c< 0$ (Theorem \ref{c<0}).
Afterwards, using Theorem \ref{c<0} and the Fredholm alternative, we can prove the existence and uniqueness of the solution to our initial problem in  the class of functions belonging to $\sol$ and having null average  (Theorem  \ref{esistenza}). 
\par Concerning the  estimate  we are led to estimate the quantity  $u-\frac{1}{|\Omega|}\int u$ instead of $u$, because the solution to the Neumann problem, if there exists, is unique up to a constant. We give three proofs of this result (Theorem \ref{stima}).  
\par We finally obtain, for the Neumann problem, similar results to those for the Dirichlet and oblique derivative boundary conditions.
\par As already said, the theorem is well known to all specialists of elliptic partial differential equations. However, after discussion with several of them, we were not able to find a precise reference for such a result, specially the one concerning the estimate.\\
\newline
This paper is organized as  follows. In Section \ref{GD} we remind the main definitions  used in this work. In Section \ref{Nex} we  prove  an existence theorem  for the Neumann problem for the Poisson's equation (Theorem \ref{esistenza}). In Section \ref{Nes} we  prove  the  estimate (Theorem \ref{stima}). In Section \ref{appendix} we remind some useful results to prove the existence theorem. 
\par We adopt the same notation used in \cite{GT}. We refer to \cite{GN} for a more detailed analysis of the problem.

\section*{Acknowledgments}

We thank  A.Adimurthi,  N.Fusco, R.Gianni, L.Orsina, and  N.Trudinger, for their advice.
We  thank B. Dacorogna who supervised the researches on this subject. 
\section{Definitions and notation}\label{GD}

In the following we denote by  $\Omega$ an open bounded non-empty set of $\mathbb{R}^N$ ($N\geq 2$)
and let $u$ be a function  defined on $\Omega$.
For every multi-index $\beta=(\beta_{1},\dots,\beta_{N})$ ($\beta_{i}\geq 0$ for $i=1,\dots,N$)  of length  $|\beta|=\sum_{i=1}^{N}\beta_{i}$, we set 
$$\displaystyle{D^{\beta}u}=\displaystyle{\frac{\partial^{|\beta|}u}{\partial x_{1}^{\beta_{1}}\dots\partial x_{N}^{\beta_{N}}}}.$$ 
We remind the definition of the usual functional spaces ($k\geq 0$):
$$C^{k}(\Omega)=\{u: \Omega \rightarrow \mathbb{R}|\mbox{\;}\forall \beta \mbox{ multi-index, }|\beta|\leq k,\,D^{\beta}u\mbox{ is continuous in } \Omega\},$$
$$C^{k}(\overline{\Omega})=\{u\in C^{k}(\Omega)| D^\beta u, \,|\beta|\leq k,\,\mbox{\;can be extended by continuity to\;} \front\}.$$
Moreover  $C^{k}(\overline{\Omega})$ is a Banach space equipped with the norm
$$\|u\|_{C^{k}}=\sum_{i=0}^{k}\sup_{|\beta|=i}\sup_{\overline{\Omega}}|D^{\beta}u|.$$

We now remind the definition of H\"older spaces.

\begin{defi}[{\bf H\"older-continuity}]
  Let $\alpha\in]0,1]$ and $\Omega\subset \mathbb{R}^N$  an open set. We define the $\alpha$-H\"older coefficient of $u: \Omega\rightarrow \mathbb{R}$ as
  $$[u]_{0,\alpha; \Omega}=\sup_{\stackrel{x,y\in\Omega}{x\neq y}}\frac{|u(x)-u(y)|}{|x-y|^{\alpha}}.$$  
If $[u]_{0,\alpha;\Omega}<\infty$ then we call u  H\"older continuous with exponent  $\alpha$ in $\Omega$.
If there is not ambiguity about the domain $\Omega$ we denote $[u]_{0,\alpha; \Omega}$ by $[u]_{0,\alpha}$.
\end{defi}

We define the {\em H\"older space} \hold\, as the set of functions belonging to  \cont\, whose $k$th-order partial derivatives are H\"older continuous with exponent  $\alpha$ in \ins.
\par  \hold \,is a Banach space equipped with the following norm:
$$\|u\|_{C^{k,\alpha}}=\|u\|_{C^{k}}+[u]_{k,\alpha}$$
where
$$[u]_{k,\alpha}=\sup_{|\beta|=k}[D^{\beta}u]_{0,\alpha}.$$

We set $C^{k,0}(\overline{\Omega})=$\cont, and one can easily verify that $C^{k,\alpha}(\overline{\Omega})\subset C^{h,\alpha}(\overline{\Omega})$ for every $h,k$ integer with $h<k$.

Moreover, for every vector-valued function $u:\Omega\rightarrow \mathbb{R}^N$, we say that it belongs to $C^{k,\alpha}(\overline{\Omega})$ if its components belong to $C^{k,\alpha}(\overline{\Omega})$.

In order to give the definition of  $C^{k,\alpha}$-domain, we need the following notations:
$$B=\{x=(x_{1},\dots,x_{N}) \in \mathbb{R}^{N}||x|<1\}\,,$$
$$B_{+}=\{x\in B|x_{N}>0\}\,, \quad B_{0}=\{x\in B|x_{N}=0\}.$$

\begin{defi}\label{Holder-domain}
We call domain every  open, bounded, connected, and non-empty subset $\Omega$ of \rn.
Moreover,  \ins\, is said to be a $C^{k,\alpha}$-domain ($k\geq 1$, $\alpha \in [0,1]$) if   for every  $p\in \partial \Omega$ there exists a neighborhood $U_p$ of $p$ in $\mathbb{R}^{N}$ and a diffeomorphism
$\varphi_p: B\rightarrow U_p$ such that
\begin{itemize}
  \item[(i)] $\varphi_p \in C^{k,\alpha}(\overline{B})\,\,and \,\,\varphi_p^{-1}\in C^{k,\alpha}(\overline{U_p})$;
  \item[(ii)] $\varphi_p(B_{+})= U_p\cap \Omega $;
  \item[(iii)] $\varphi_p(B_{0}) = U_p\cap \partial{\Omega}$.
\end{itemize}
\end{defi}

\begin{os}[{\bf H\"older-continuity on the boundary}] Let $\Omega$ be a $C^{k,\alpha}$-domain.
We say that $u\in C^{k,\alpha}(\partial \Omega)$ if, for every $p\in \partial \Omega$, we have $u\circ \varphi_p \in C^{k,\alpha}(B_0)$, where  $\varphi_p$ is given by the previous definition. 

Of course if $u\in C^{k,\alpha}(\overline{\Omega})$ its restriction to $\partial \Omega$ belongs to $C^{k,\alpha}(\partial \Omega)$. 
Moreover, for every $u\in C^{k,\alpha}(\partial \Omega$) there exists a function belonging to  $C^{k,\alpha}(\overline{\Omega})$ whose restriction to the boundary coincides with $u$ (see \cite{GT}: Lemma 6.38 p.137). 

Then, in order to study the Neumann problem for the Poisson's equation, we can consider boundary values belonging to  $C^{k,\alpha}(\overline{\Omega})$ instead of $ C^{k,\alpha}(\partial \Omega)$.

\end{os}


We denote by $L^p$ and $W^{m,p}$ the usual Lebesgue and Sobolev spaces and we refer to \cite{A, B}  for their properties.
\section{Existence of solutions to the Neumann problem for the Poisson's equation in $\sol$}\label{Nex}

In this section we consider the Poisson's equation with Neumann boundary condition and we  prove the following result:

\begin{teo}\label{esistenza}
Let  $\Omega$ be a  $C^{2,\alpha}$-domain and let  $f\in \eqz,g\in \borneu$ such that:
\begin{equation}\label{comp}
\int_{\Omega} f=\int_{\front} g.
\end{equation}
Then the problem
\begin{equation}\label{problem}
\left\{
\begin{array}{rl}
\Delta u=f &\textrm{ in $\Omega$} \\
\dn=g &\textrm{ on $\front$}
\end{array} \right.
\end{equation}
admits a unique solution in the class
$$\mathscr{C}=\left\{ u\in \sol:\frac{1}{|\Omega|}\int_{\Omega} u=0 \right\}.$$
\end{teo}

We start by reminding the following estimate which is very usefull in the following:

\begin{teo}\label{stima Fiorenza}
Let \ins\, be a $C^{2,\alpha}$-domain and let $u\in \sol$ be a solution of \eqref{problem}. Then
\begin{equation}
\label{stima Neumann}
\|u\|_{C^{2,\alpha}}\leq C(\|f\|_{C^{0,\alpha}}+\|g\|_{C^{1,\alpha}}+\|u\|_{C^{0}}),
\end{equation}
with $C=C(\Omega,\alpha,N).$\\
(\cite{GT}: Theorem 6.30, p.127;  \cite{LU}: Theorem 3.1, p.126).
\end{teo}

We prove two preliminary lemmas.

\begin{lemma}
Let \ins$\subset \mathbb{R}^{N}$ be a $C^{2}$-domain and  let $u\in C^{2}(\overline{\Omega})$. We suppose that there exists  $p \in \front$ such that:
$$u(p)=\max_{\overline{\Omega}}u\mbox{\;}(\min_{\overline{\Omega}}u)\mbox{\;\;\;and\;\;\;}Du(p)=0.$$
then  $\Delta u(p)\leq 0\mbox{\;}(\geq 0)$.
\end{lemma}

\dimo:  Up to a translation, we can suppose that $p$ is the origin of \rn. Moreover, as $\Omega \in C^{2}$, up to a rotation, we can also suppose that, for some $R$, we have $-re_i\in \Omega$ for every $r\in [0,R]$ and $i=1,...,N$ where $(e_1,...,e_N)$ is an orthonormal basis of  \rn.
Then, partial derivatives  verify     
$$D_{ii}u(p)=\lim_{h\rightarrow 0^{-}}\frac{D_{i}u(p+he_{i})-D_{i}u(p)}{h}=\lim_{h\rightarrow 0^{-}}\frac{D_{i}u(p+he_{i})}{h}\leq 0,$$
because $Du(p)=0$ and, for $h<0$ small enough we have $D_{i}u(p+he_{i})\geq 0$.\\
If $p$ is a minimum point the proof is similar.  
$\Box$

\begin{lemma}
Let $\Omega\subset \mathbb{R}^{N}$ be a $C^{2}$-domain and  $f\in \eqz$.
Let $u\in \sol$ be a solution to the problem 
$$
\left\{
\begin{array}{rl}
\Delta u -u=f &\textrm{ in $\Omega$} \\
\dn=0 &\textrm{ on $\front$}
\end{array} \right.
$$
then
$$\|u\|_{C^{0}}\leq \|f\|_{C^{0}}.$$
\end{lemma}

\dimo: Let $p$ be such that  $|u(p)|=\max_{\overline{\Omega}}|u|$.  We consider separately the following cases: $p\in \Omega$ and $p\in \front$.\\
${\bf p\in \Omega}$. Let $p$ be a maximum point  for   $u$ then $u(p)\geq 0$  and  $\Delta u(p)\leq 0$. As $\Delta u -u =f$ we get $f(p)=\Delta u(p)-u(p)\leq \Delta u(p)\leq 0$.
Thus
$$\|u\|_{C^{0}}=u(p)=\Delta u(p)-f(p)\leq -f(p) \leq \|f\|_{C^{0}}.$$
If  $p$ is a minimum point for $u$ ($u(p)\leq 0$) the proof is similar. \\
${\bf p\in \front}$. If $p$ is a maximum point for  $u$ then $p$ is a maximum point for  $u_{|\front}$. So, for every tangent vector  $\tau$ to $\front$ at $p$ we have:
$$\frac{\partial u}{\partial \tau}(p)=0$$
and, using the hypothesis 
$$\dn(p)=0,$$
we get $Du(p)=0$. Thanks to the previous lemma we have  $\Delta u(p)\leq 0$ and, as $\Delta u-u=f$ on $\front$, the lemma ensues.
We can use the same arguments if  $p$ is a  minimum point for $u$.
$\Box$

Now, we can show the existence theorem for the Neumann problem for the Laplace operator:

{\bf Proof of Theorem \ref{esistenza}}: By Theorem \ref{c<0}  there exists an unique solution  in  $\sol$ , denoted by $\textfrak{T}[f,g]$, to the problem:
$$
\left\{
\begin{array}{rl}
\Delta u-u=f &\textrm{ in $\Omega$} \\
\dn=g &\textrm{ on $\front$}
\end{array} \right.
$$
for every $f,g$ verifying the compatibility condition \eqref{comp}. Moreover 
\begin{equation}\label{ok}
\int_{\Omega} u=\int_{\Omega} \Delta u-\int_{\Omega}f=\int_{\front} \displaystyle{\frac{\partial u}{\partial n}}-\int_{\front} g=0\,,
\end{equation}
so $\textfrak{T}[f,g] \in \mathscr{C}$. Defining 
$$\mathscr{A}=\left\{(f,g)\in \eqz \times \borneu : \int_{\Omega}f=\int_{\partial \Omega }g\right\}$$ 
we have that   $\textfrak{T}: \mathscr{A} \rightarrow \mathscr{C}$ is a well defined bijective operator.
\par Now, we consider the following equation:
\begin{equation}\label{eqz}
u- \textfrak{T}[-u,0]= \textfrak{T}[f,g].
\end{equation}
\indent  Then  $u \in \mathscr{C}$ is a solution to the problem  \eqref{problem} if and only if it is a solution of \eqref{eqz}, because:
$$(f,g)=\textfrak{T}^{-1}(u- \textfrak{T}[-u,0])=\textfrak{T}^{-1}u+(u,0)=\left(\Delta u,\dn\right).$$ 
Then, we need to show that for every  $(f,g)\in \mathscr{A}$ there exists an unique solution $u\in \mathscr{C}$ of \eqref{eqz}. As $\textfrak{T}$ is bijective we are led to prove that, for every   $v\in  \mathscr{C}$, the equation
\begin{equation}\label{eqz1}
u-\textfrak{T}[-u,0]=v
\end{equation} 
admits an unique solution on  $\mathscr{C}$. For that we use the Fredholm alternative theorem.\\
\indent We consider the space 
$$\mathscr{F}=\left\{f\in \eqz : \frac{1}{|\Omega|}\int_{\Omega} f=0\right\},$$
equipped with the norm of $\eqz$.
Let $T$ be the operator:
$$T: \mathscr{F}\rightarrow \mathscr{F},$$
$$Tf=\textfrak{T}[-f,0].$$
Using the properties of  $\textfrak{T}$, we get: 
\begin{equation}\label{eqz0}
T(\mathscr{F}) \subset \mathscr{C}.
\end{equation}
\indent We firstly show that  $T$ is a compact operator. 
Let $\{f_{k}\}\subset \mathscr{F}$ then, because of Theorem \ref{c<0} and \eqref{ok}, there exists $\{u_{k}\}\subset \mathscr{C}$ such that $Tf_{k}=u_{k}$ and, because of \eqref{stima Neumann}, we have
$$\|u_{k}\|_{C^{2,\alpha}}\leq C(\Omega, \alpha, N) (\|f_{k}\|_{C^{0,\alpha}}+\|u_{k}\|_{C^{0}}).$$
So, using the previous lemma, we get:
$$\|u_{k}\|_{C^{2,\alpha}}\leq 2C(\Omega, \alpha, N)\|f_{k}\|_{C^{0,\alpha}}.$$
If $\{f_{k}\}$ is a bounded sequence  of $\mathscr{F}$ then $\{u_{k}\}$ is  bounded in $\sol$ and in $W^{2,\infty}(\Omega).$
Thus there exists a subsequence  $\{u_{k_{h}}\}$ and a function  $u\in W^{2,\infty}(\Omega)$ such that
$$u_{k_{h}}\stackrel{*}{\rightharpoonup} u,\mbox{\;\;in\;\;}W^{2,\infty}$$
so, for every $p>1$
$$u_{k_{h}}\stackrel{*}{\rightharpoonup} u,\mbox{\;\;in\;\;}W^{2,p}.$$
Now, choosing $\displaystyle{\frac{N}{2-\alpha}}<p\leq N$, by Rellich-Kondrachov theorem, we obtain
$$u_{k_{h}}\rightarrow u \mbox{\,\,\, in\,\,\, } C^{0,\alpha}(\overline{\Omega})$$
which proves that $T$ is compact.\\
\indent Now, equation \eqref{eqz1} can be rewritten as
\begin{equation}\label{eqz2}
u-Tu=v.
\end{equation}
The equation
$$u-Tu=0$$
is equivalent to the problem
$$
\left\{
\begin{array}{rl}
\Delta u=0 &\textrm{ in $\Omega$} \\
\dn=0 &\textrm{ su $\front$}
\end{array} \right.
$$
which, in $\mathscr{C}$, admits only the trivial solution  $u=0$.
\par Then, as $T$ is compact, applying to \eqref{eqz2} the Fredholm alternative we have that for every  $v\in \mathscr{C}$ there exists an unique solution  $u\in \mathscr{C}$ of \eqref{eqz2}  and the theorem ensues.
$\Box$

\begin{os}
If  $u\in \mathscr{C}$ is a solution to the problem  \eqref{problem} then for every  $k \in \mathbb{R}$ the function $u+k$ is also a solution to  \eqref{problem} in $\sol$ (but not in $\mathscr{C}$). \\
On the other hand, if  $u$ is a solution to  \eqref{problem} in $\sol$ we can obtain a solution in  $\mathscr{C}$ setting
$$v=u-\frac{1}{|\Omega|}\int_{\Omega} u.$$
Then, using Theorem \ref{esistenza} we have a existence and uniqueness (up to an additive  constant) result for the problem \eqref{problem} in $\sol$.
\end{os}
\section{Schauder estimate}\label{Nes}

We prove an estimate of  $\|u\|_{C^{2,\alpha}}$ in terms of  $\|f\|_{C^{0,\alpha}},\|g\|_{C^{1,\alpha}}$ for any  solution to problem \eqref{problem}. In particular, as the uniqueness of the solution is proved  up to an additive constant,  we prove the  estimate for a solution with null average. \\
\indent We state the following theorem and we give three proofs:
\begin{teo}\label{stima}
Let $\Omega$ be a $C^{2,\alpha}$-domain and   $f\in C^{0,\alpha}(\overline{ \Omega})$, $g\in\borneu$ such that
$$\int_{\Omega}f=\int_{\partial \Omega}g.$$ Let $u\in \sol$ be a solution to \eqref{problem}.
Then
$$\left\|u-\frac{1}{|\Omega|}\int_{\Omega} u\right\|_{C^{2,\alpha}}\leq C(\|f\|_{C^{0,\alpha}}+\|g\|_{C^{1,\alpha}}),$$
with $C=C(\Omega,\alpha, N)$.
\end{teo}

The first proof of Theorem \ref{stima} has been suggested to us by A. Adimurthy who attributed it to Jacques-Louis Lions. However we did not find any reference for such a proof, which is detailed in the following.

{\bf First proof of Theorem \ref{stima}}: Let  $u\in \sol$ be a solution to  \eqref{problem} such that $\frac{1}{|\Omega|}\int_{\Omega} u=0$ . By  \eqref{stima Neumann}, we have the following estimate for $u$ 
\begin{equation}\label{vecchia}
\|u\|_{C^{2,\alpha}}\leq C_{1}(\Omega, \alpha,N)(\|u\|_{C^{0}}+\|f\|_{C^{0,\alpha}}+\|g\|_{C^{1,\alpha}})
\end{equation}
and we are interested in proving that:
\begin{equation}\label{nuova}
\|u\|_{C^{2,\alpha}}\leq C_{2}(\Omega, \alpha,N)(\|f\|_{C^{0,\alpha}}+\|g\|_{C^{1,\alpha}}).
\end{equation}
\indent Let us suppose that  \eqref{nuova} is false. Then for every $k\in\mathbb{N}$ there exists  $\{u_{k}\}\in \sol$ and  $\{f_{k}\}\in\eqz,\{g_{k}\}\in\borneu$ such that
\begin{equation}\label{1}
\left\{
\begin{array}{rl}
\Delta u_{k}=f_{k} &\textrm{ in $\Omega$} \\
\displaystyle{\frac{\partial u_{k}}{\partial n}}=g_{k} &\textrm{ on $\front$}
\end{array} \right.
\end{equation}
\begin{equation}\label{2}
\frac{1}{|\Omega|}\int_{\Omega} u_{k}=0\,,
\end{equation}
\begin{equation}\label{3}
\|u_{k}\|_{C^{2,\alpha}}=1\,,
\end{equation}
\begin{equation}\label{4}
\|u_{k}\|_{C^{2,\alpha}}>k(\|f_{k}\|_{C^{0,\alpha}}+\|g_{k}\|_{C^{1,\alpha}})\,.
\end{equation}
Thus, we have   $f_{k}\rightarrow 0$ in $\eqz$ and $g_{k}\rightarrow 0$ in $\borneu$.\\
\indent Using \eqref{3} we have that for every multi-index $\beta$, $|\beta|=0,1,2$, $\{D^{\beta}u_{k}\}$ is uniformly bounded in $C^{0}(\overline{\Omega})$ and equicontinuous because 
$$|D^{\beta}u_{k}(x)-D^{\beta}u_{k}(y)|\leq |x-y|^{\alpha}\mbox{\;\;\;}\forall x,y\in\Omega\mbox{, }\forall |\beta|=2$$
which implies that
$$|D^{\beta}u_{k}(x)-D^{\beta}u_{k}(y)|\leq C(\Omega)|x-y|^{\alpha}\mbox{\;\;\;}\forall x,y\in\Omega\mbox{, }\forall |\beta|=0,1.$$
Iterating the Ascoli-Arzel\`a theorem  we get a subsequence $\{u_{k_{h}}\}$ such that
\begin{equation}
\begin{array}{l}
u_{k_{h}}\rightarrow u_{0}\mbox{\;\;in\;\;}C^{0}(\overline{\Omega}),\vspace{0.2cm}\\

D^{\beta}u_{k_{h}}\rightarrow u_{\beta}\mbox{\;\;in\;\;}C^{0}(\overline{\Omega})\mbox{\;\;}\forall \beta, |\beta|=1,2,
\end{array}
\end{equation}
which implies
$$u_{k_{h}}\rightarrow u_{0}\mbox{\;\;in\;\,}C^{2}(\overline{\Omega}).$$
Then
$$\Delta u_0=\lim_{h}\Delta u_{k_{h}}=\lim_h f_{k_{h}}=0$$
$$\displaystyle{\frac{\partial u_0}{\partial n}=\lim_h \frac{\partial u_{k_{h}}}{\partial n}=\lim_h g_{k_{h}}=0}$$
and
$$\left\{
\begin{array}{rl}
\Delta u_{0}=0 &\textrm{ in $\Omega$} \\
\displaystyle{\frac{\partial u_{0}}{\partial n}}=0&\textrm{ on $\front$}\\
\frac{1}{|\Omega|}\int_{\Omega} u_{0}=0 &    \\
\end{array} \right.
$$
which implies  $u_{0}= 0$. Comparing with \eqref{vecchia}, we get a contradiction because  
$$1=\|u_{k_{h}}\|_{C^{2,\alpha}}\leq C_{1}(\Omega, \alpha,N)(\|u_{k_{h}}\|_{C^{0}}+\|f_{k_{h}}\|_{C^{0,\alpha}}+\|g_{k_{h}}\|_{C^{1,\alpha}})\rightarrow 0.$$
$\Box$ 

For the second proof of  Theorem \ref{stima} we need two more lemmas. The first one states an estimate  in $L^2(\Omega)$ for $u-\frac{1}{|\Omega|}\int_{\Omega} u$:
\begin{lemma}\label{debole}
Let $\Omega$ be a $C^{2,\alpha}$-domain and $f\in C^{0,\alpha}(\overline{ \Omega})$, $g\in\borneu$ such that
$$\int_{\Omega}f=\int_{\partial \Omega}g.$$
Let $u\in \sol$ be a  solution to  \eqref{problem}.
Then
$$\left\|u-\frac{1}{|\Omega|}\int_{\Omega} u\right\|_{L^{2}}\leq C(\|f\|_{C^{0,\alpha}}+\|g\|_{C^{1,\alpha}}),$$
with $C=C(\Omega,N)$.
\end{lemma}
\dimo: We can suppose  $\frac{1}{|\Omega|}\int_{\Omega} u=0$. Integrating by parts we have  
$$\int_{\Omega}|Du|^{2}=\int_{\partial \Omega}ug-\int_{\Omega}uf$$
and using the Young inequality   with $\varepsilon>0$  we get 
$$\int_{\Omega}|Du|^{2}\leq \int_{\partial \Omega}|ug|+\int_{\Omega}|uf|\leq\varepsilon\int_{\partial \Omega}|u|^{2}+\frac{1}{4\varepsilon}\int_{\partial \Omega}|g|^{2}+\varepsilon\int_{\Omega}|u|^{2}+\frac{1}{4\varepsilon}\int_{\Omega}|f|^{2}.$$
Now,  we have
$$\int_{\partial \Omega}|u|^{2}\leq C_{1}(\Omega)\|u\|_{W^{1,2}}^{2}\mbox{\;\;}\forall u\in W^{1,2}(\Omega)\cap C(\overline{\Omega});$$
so
$$\int_{\Omega}|Du|^{2}\leq \varepsilon C_{2}(\Omega)\|u\|_{W^{1,2}}^{2}+C_{3}(\Omega,\varepsilon)\left[\|g\|_{W^{1,2}}^{2}+\|f\|_{L^{2}}^{2}\right]$$
and, as $u\in \mathscr{C}$, using the Poincar\' e Inequality, we get
$$\int_{\Omega}|Du|^{2}\leq \varepsilon C_{4}(\Omega,N)\|Du\|_{L^{2}}^{2}+C_{5}(\Omega,\varepsilon)\left[\|g\|_{C^{1,\alpha}}^{2}+\|f\|_{C^{0,\alpha}}^{2}\right].$$
Choosing $\varepsilon<1/C_{4}(\Omega,N)$ and using the Poincar\'e Inequality we have
$$\|u\|_{L^{2}}\leq C_{6}(\Omega,N)\|Du\|_{L^{2}}\leq C\left(\|g\|_{C^{1,\alpha}}+\|f\|_{C^{0,\alpha}}\right),$$
with $C=C(\Omega,N)$.
$\Box$

The second lemma proves a local estimate for solution of Poisson's equation:
\begin{lemma}\label{locale}
Let $\Omega$ be a domain of \rn. Let  $f\in L^{p}(\Omega)$ with $p>N/2$ and  $u\in C^{2}(\Omega)$ a solution of $\Delta u=f$ in \ins.  Then, for every ball
$B(y,2R)\subset \Omega$ we have 
$$\sup_{B(y,R)}|u|\leq C\left(R^{-\frac{N}{2}}\|u\|_{L^{2}(B(y,2R))}+R^{2-\frac{N}{p}}\|f\|_{L^{p}(\Omega)}\right),$$
with $C=C(N,p)$.\\
(\cite{SE}: Theorem 1 p.255 and  Theorem 2 p.259).
\end{lemma}

Following a N. Fusco's idea, we can now give the second proof of Theorem \ref{stima}:\\
{\bf Second proof of Theorem \ref{stima}}: We can suppose $\frac{1}{|\Omega|}\int_{\Omega} u=0$. We give an estimate for   $\sup_{\partial \Omega}|u|$. As  $\Omega$ is a $C^{2,\alpha}$-domain, for every   $\varepsilon$ small enough we have 
$$(x-\varepsilon n(x))\in \Omega\quad\forall x\in\front$$
where $n(x)$ denotes the unit outer-pointing  normal at $x$ to $\partial \Omega$. Then, by  Lagrange theorem, 
for every $x\in \front$ there  exists $\tau=\tau(x,\varepsilon)\in(0,1)$ such that  
$$ u(x-\varepsilon n(x))=u(x)+\left(Du(x-\tau \varepsilon n(x)),-\varepsilon n(x)\right).$$
Defining  $\Omega_{\varepsilon}=\{x\in \Omega|\mbox{\;dist}(x,\front)\geq \varepsilon\}$ and taking the supremum for $x\in \front$ in the previous relationship we get:\begin{equation}\label{frontiera}
\sup _{\front}|u|\leq \sup_{\Omega_{\varepsilon}}|u|+\varepsilon\sup_{\Omega}|Du|.
\end{equation}
Let  $R<\varepsilon/2$. We consider a finite cover of  $\Omega_{\varepsilon}$, denoted by  $\{B(y_{i},R)\}_{i=1}^{M}$ with $M=M(\varepsilon,\Omega$) ($B(y_{i},R)$ denotes an open ball  of radius  $R$ centered at $y_i$). As  $f\in C^{0}(\overline{\Omega})$ then  $f\in L^{N+1}(\Omega)$ and, by Lemma \ref{locale}, we have 
$$\sup_{\Omega_{\varepsilon}}|u|\leq \sum_{i=1}^{M}\sup_{B(y_{i},R)}|u|\leq MC_{1}(N)\left(R^{-\frac{N}{2}}\|u\|_{L^{2}(B(y,2R))}+R^{2-\frac{N}{N+1}}\|f\|_{L^{N+1}(\Omega)}\right)$$
which implies
\begin{equation}\label{compatto}
\sup_{\Omega_{\varepsilon}}|u|\leq C_{2}(\varepsilon,\Omega,N)\left(\|u\|_{L^{2}}+\|f\|_{L^{N+1}}\right).
\end{equation}
By Theorem 3.7 in \cite{GT} (p.36) and  \eqref{frontiera}, we have 
$$\sup_{\Omega}|u| \leq \sup_{\Omega_{\varepsilon}}|u|+\varepsilon\sup_{\Omega}|Du| + C_{3}(\Omega)\sup_{\Omega}|f|$$
and, by \eqref{compatto}, we get 
$$\sup_{\Omega}|u| \leq C_{2}(\varepsilon,\Omega,N)\left(\|u\|_{L^{2}}+\|f\|_{L^{N+1}}\right)+\varepsilon\sup_{\Omega}|Du| + C_{3}(\Omega)\sup_{\Omega}|f|.$$
Using Lemma \ref{debole}, we get
$$\sup_{\Omega}|u| \leq C_{4}(\varepsilon,\Omega,N)(\|f\|_{C^{0,\alpha}}+\|g\|_{C^{1,\alpha}})+\varepsilon \|u\|_{C^{2,\alpha}}$$
and, by \eqref{stima Neumann}, we have 
$$\|u\|_{C^{2,\alpha}}\leq C_{5}(\Omega,\alpha,N,\varepsilon)(\|f\|_{C^{0,\alpha}}+\|g\|_{C^{1,\alpha}})+\varepsilon C_{6}(\Omega,\alpha, N)\|u\|_{C^{2,\alpha}}.$$
Choosing  $\varepsilon<1/C_{6}(N,\Omega,\alpha)$ the theorem ensues.
$\Box$

We give now the third proof of Theorem  \ref{stima}. We thank R. Gianni who gave us  some ideas for this proof.

{\bf Third proof of Theorem \ref{stima}}: We  suppose $\int_{\Omega} u=0$ and, because of  \eqref{stima Neumann}, we just need to prove an estimate for  $\|u\|_{C^{0}}$ in terms of  $\|f\|_{C^{0,\alpha}}$ and $\|g\|_{C^{1,\alpha}}$. \\
\indent Let $M=\underset{\overline{\Omega}}{\max}|u|$ and  $p\in\overline{\Omega}$ such that $|u(p)|=M$. Moreover, we suppose  $M>(\|f\|_{C^{0,\alpha}}+\|g\|_{C^{1,\alpha}})$, otherwise there is nothing to prove. By  \eqref{stima Neumann}, we have 
\begin{equation}\label{Du}
\|Du\|_{C^{0}}\leq C_1(\|f\|_{C^{0,\alpha}}+\|g\|_{C^{1,\alpha}}+\|u\|_{C^{0}})=K\,,
\end{equation}
with $C_1= C_{1}(\Omega,\alpha,N)$. We distinguish two cases: $p\in \Omega$ and $p\in \partial \Omega$. 

In the following we consider two  cases: $p\in \Omega$ and $p\in \partial \Omega$. We suppose that  $p\in \Omega$ and we denote by $B(p,r)$ the open ball of radius $r$ centered at $p$. We consider in particular $r<r_0=\rm{dist}(p,\partial\Omega)$ so that  such a ball is contained in $\Omega$.


 Now, we can prove that, for every $r<r_0$, we have 
$$|u|>M-rK\mbox{,\quad in\;\;}
B(p,r).$$
In fact, for every $x\in B(p,r)$ we have 
$$\displaystyle{u(x)-u(p)=\int_{0}^{1}\frac{d}{dt}u(p+t(x-p))\, dt=\int_{0}^{1}\langle Du(p+t(x-p)),x-p\rangle\,dt}$$
so
$$|u(x)|\geq |u(p)|-Kr=M-Kr.$$
\indent Choosing  $r<{\rm min}\{1/4C_{1},r_0\}$, for every  $x\in B(p,r)$, we have:
$$
\begin{array}{ll}
|u(x)|\geq M-rK &=(1-rC_{1})M-rC_{1}(\|f\|_{C^{0,\alpha}}+\|g\|_{C^{1,\alpha}})>\\
& \displaystyle{>\frac{3}{4}M-rC_{1}(\|f\|_{C^{0,\alpha}}+\|g\|_{C^{1,\alpha}})}.\\
\end{array}
$$
Then, as  $M>(\|f\|_{C^{0,\alpha}}+\|g\|_{C^{1,\alpha}})$ and $r<1/4C_1$, we get 
$$|u(x)|>\frac{3}{4}M-C_{1}rM>\frac{M}{2}\mbox{\;\;\;}\forall x\in B(p,r)$$
and, denoting by $\omega_N$ the measure of the unit ball of \rn, we get
$$\frac{M}{2}\frac{\omega_{N}^{1/2}}{2^{N/2}}r^{N/2}\leq \left(\int_{B(p,r)}|u|^{2}\right)^{1/2}\leq \|u\|_{L^{2}}.$$
\indent Then, by Lemma \ref{debole},  we can state
$$\|u\|_{C^{0}}\leq \frac{2^{\frac{N+2}{2}}r^{-\frac{N}{2}}}{\omega_{N}^{1/2}}C_{2}(\Omega,N)(\|f\|_{C^{0,\alpha}}+\|g\|_{C^{1,\alpha}})$$
which implies
$$\|u\|_{C^{2,\alpha}}\leq C(\|f\|_{C^{0,\alpha}}+\|g\|_{C^{1,\alpha}}),$$
with $C=C(\Omega,\alpha,N).$

We suppose now that  $p\in \partial \Omega$. As $\Omega$ is a $C^2$-domain, it satisfies  the interior sphere condition (see \cite{GT} p.33). Then, there exists $r_0>0$ such that the open ball $B(p-r n(p),r)$ with $r<r_0$  is contained in $\Omega$  ($n(p)$ denotes the unit outer-pointing normal at $p$ to the boundary). 

Thus, as $p$ is not the center of the ball, similarly to the previous case, we obtain
$$|u|>M-2rK\mbox{\quad in\;\;}
B(p-rn(p),r)\quad r<r_0\,.$$
Then, by choosing $r<{\rm min}\{1/8C_{1},r_0\}$, 
we get the result by the same argument used above applied to the ball $B(p-r n(p),r)$. 
$\Box$

\section{Appendix: a preliminary result}\label{appendix}

In this section we remind some result used in Section \ref{Nex}. We consider a $C^{2,\alpha}$-domain $\Omega$ and the following operator 
$$L:C^{2, \alpha}(\overline{\Omega})\rightarrow \mathbb{R} $$
$$Lu=\sum_{i,j=1}^{N}a_{ij}(x)\frac{\partial^{2}{u}}{\partial{x_{i}}\partial{x_{j}}}+\sum_{i=1}^{N}b_{i}(x)\frac{\partial u}{\partial x_{i}}+c(x)u$$
where $a_{ij}, b_i, c\in C^{0,\alpha}(\overline{\Omega})$ for every $i,j=1,...,N$  and  $c(x)\leq 0$ for every $x\in \overline{\Omega}$.

We say that $L$ is elliptic is the matrix $A(x)=[a_{i,j}(x)]$ is positively defined for every $x\in \overline{\Omega}$. Moreover, denoting $\lambda(x)$  the smallest eigenvalue of $A(x)$, we said that $L$ is uniformly elliptic if there exists $\lambda_0>0$ such that 
$\lambda(x)\geq \lambda_0$ for every $x\in \overline{\Omega}$.

We refer to  \cite{GT} (p. 130) for the following alternative result for the oblique derivative problem:

\begin{teo}\label{alternativa}
Let $\Omega$ be a $C^{2,\alpha}$-domain. Let  $L$ be an elliptic operator, which is uniformly elliptic in $\Omega$ and with $C^{0,\alpha}(\overline{\Omega})$-coefficients. Let   $l,m\in \borneu$ such that $m \neq 0$ for every $x\in \partial{\Omega}$.
Then exactly one of the following holds:

\begin{itemize}
 \item[1.] the homogeneous problem
$$
\left\{
\begin{array}{rl}
Lu=0 &\textrm{ in $\Omega$} \\
l(x)u+m(x)\dn=0 &\textrm{ on $\front$}
\end{array} \right.
$$
admits nontrivial  solutions;
 \item[2.] homogeneous problem has only the trivial solution, in which case for every $f\in \eqz, g\in \borneu$ there exists  a solution  $u\in \sol$ to the inhomogeneous  problem
$$
\left\{
\begin{array}{rl}
Lu=f &\textrm{ in $\Omega$} \\
l(x)u+m(x)\dn=g &\textrm{ on $\front$}
\end{array} \right.
$$ 
\end{itemize}

\end{teo}

We can now prove the following theorem:

\begin{teo}\label{c<0}
Let  \ins\, be $C^{2,\alpha}$-domain. Then for every $f\in \eqz, g\in \borneu$ there exists an unique solution $u\in \sol$ to the problem
$$
\left\{
\begin{array}{rl}
\Delta u -u=f &\textrm{ in $\Omega$} \\
\dn=g &\textrm{ on $\front$}
\end{array} \right.
$$
\end{teo}
\dimo: The problem
$$
\left\{
\begin{array}{rl}
\Delta u -u=0 &\textrm{ in $\Omega$} \\
\dn=0 &\textrm{ on $\front$}
\end{array} \right.
$$
admits only the trivial solution. Using point 2 of Theorem  \ref{alternativa} the result ensues.
$\Box$

\bibliographystyle{plain}  
\bibliography{bibliography}	

\end{document}